\documentclass[11pt]{amsart}
\usepackage{geometry,amssymb}
\usepackage{hyperref}
\usepackage{mathtools}
\newcommand{\IZ}{\mathbb Z}

\theoremstyle{plain}
\newtheorem{proposition}{Proposition}
\newtheorem{theorem}{Theorem}

\theoremstyle{definition}

\begin{document}
    \title{Steiner systems $S(2,6,226)$ and $S(2,6,441)$ exist}
    \author{Taras Banakh, Ivan Hetman, Alex Ravsky}
    \address{T.~Banakh: Ivan Franko National University of Lviv (Ukraine),\newline and Jan Kochanowski University in Kielce (Poland)}
    \address{I.~Hetman: Lviv (Ukraine)}
    \address{A.~Ravsky: Pidstryhach Institute for Applied Problems of Mechanics and Mathematics
    National Academy of Sciences of Ukraine, Lviv (Ukraine)}
    \email{t.o.banakh@gmail.com, vespertilion@gmail.com, vitaravski@gmail.com}

    \begin{abstract} Via computer search, we found seven non-isomorphic $1$-rotational Steiner systems $S(2,6,226)$ and six point-transitive Steiner systems $S(2,6,441)$, resolving two of $29$ previously undecided cases for $S(2,6,v)$.
    \end{abstract}

    \maketitle

    \section{Introduction}\label{sec:intro}

    A {\em Steiner system} $S(2,k,v)$ is a set $V$ of $v$ \emph{points} endowed with a family $\mathcal B$ of $k$-element subsets of $V$, called \emph{blocks}, such that every two distinct points of $V$ belong to a unique block.
    There are well-known necessary conditions of the existence of $S(2,k,v)$.

    \begin{proposition}
        If a Steiner system $S(2,k,v)$ exists, then $k-1 | v-1$, $k(k-1) | v(v-1)$, and $v\ge k^2 - k + 1$.
    \end{proposition}

    Numbers $v$ satisfying the necessary conditions above are called \emph{admissible} (with respect to the number $k$).
    By $b=\frac{v(v-1)}{k(k-1)}$ we denote the number of blocks.
    For $k\le 5$, Steiner systems $S(2,k,v)$ exist for all admissible $v$ according to the results of Reiss~\cite{Reiss}, Moore~\cite{Moore} and Hanani~\cite{Hanani}.
    For $k\ge 6$, the existence of $S(2,k,v)$ was proved~\cite{Wilson} only for sufficiently large admissible $v$.

    \begin{theorem}[Wilson, 1975]
        For every $k\ge 3$, there exists a number $v_k$ such that for every admissible number $v \ge v_k$ there exists a Steiner system $S(2,k,v)$.
    \end{theorem}

    On the other hand, the famous Bruck--Ryser Theorem~\cite{BrRys} implies that no Steiner system $S(2,6,36)$ exists, and computer search~\cite{HTJL} showed that no $S(2,6,46)$ exists.
    Denniston, Mills, Wilson, and others contributed many examples reducing the number of undecided cases.
    The Handbook of Combinatorial Designs~\cite[\S IV.3.4]{HB} lists $29$ numbers \[51, 61, 81, 166, 226, 231, 256, 261, 286, 316, 321, 346, 351, 376, 406,\] \[411, 436, 441, 471, 501, 561, 591, 616, 646, 651, 676, 771, 796, 801,\] for which the existence of $S(2,6,v)$ is not known.
    In this paper we present Steiner systems $S(2,6,226)$ and $S(2,6,441)$ thus resolving two of those 29 undecided cases.
    By computer search, we have found seven Steiner systems $S(2,6,226)$ as $1$-rotational difference families on the $225$-element group $(\IZ_5\times\IZ_5\times\IZ_3)\rtimes\IZ_3$ with an attached point, and six Steiner systems $S(2,6,441)$ as difference families on the $441$-element group $(\IZ_7\rtimes\IZ_3)\times(\IZ_7\rtimes\IZ_3)$.

    Difference families are convenient and compact way to present Steiner systems.
    Given a group $G$ of cardinality $v$, we look for a family $\mathfrak B$ of $k$-element subsets of $G$ such that the pair $(G, \{Bg:B\in\mathfrak B,g\in G\})$ is a Steiner system $S(2,k,v)$.
    In this case, the family $\mathfrak B$ is called a $(v,k,1)$ {\em difference family}, and its elements are called {\em base blocks}, see~\cite[\S VI.16.1]{HB}.

    A slight modification of this construction leads to the following notion of a $1$-rotational difference family.
    Given any group $G$ of cardinality $v-1$, choose any point $\infty\notin G$ and extend the right action of the group $G$ on itself to a right action of $G$ on $G\cup\{\infty\}$, letting $\infty{\cdot} g=\infty$ for all $g\in G$.
    Then we look for a family $\mathfrak B$ of $k$-element subsets of $G\cup\{\infty\}$ such that the pair $(G\cup\{\infty\},\{Bg:B\in\mathfrak B,g\in G\})$ is a Steiner system $S(2,k,v)$.
    In this case, the family $\mathfrak B$ is called a $(v,k,1)$ {\em $1$-rotational difference family}, and its elements are called {\em base blocks}, see~\cite[\S VI.16.6]{HB}.

    For a subset $B\subseteq G\cup \{\infty\}$, the subgroup $G_B\coloneqq\{g\in G:B=Bg\}$ is called the \emph{stabilizer} of the set $B$.
    It is well-known that $|G|=|G_B|\cdot |Orb_GB|$, where $Orb_GB\coloneqq\{Bg:g\in G\}$ is the \emph{orbit} of the set $B$ under the action of the group $G$.

    We denote by $\IZ_n$ the cyclic group $\{0,1,\dots,n-1\}$ of order $n$ under the addition modulo $n$.
    In order to present Steiner systems $S(2,6,226)$ and $S(2,6,441)$, we use the groups $(\IZ_5\times \IZ_5\times\IZ_3)\rtimes\IZ_3$ and $(\IZ_7\rtimes\IZ_3)\times(\IZ_7\rtimes \IZ_3)$ of orders $225$ and $441$, respectively.

    The algorithm used to produce the Steiner systems $S(2,6,226)$ and $S(2,6,441)$ is described in the presentation~\cite{Het} (and a subsequent paper) of the second author.

    \section{Steiner systems $S(2,6,226)$}\label{sec:226}

    The construction described below is similar to that of Mills~\cite{WHM78b} for $S(2,6,136)$.

    The group whose one-point extension contains $1$-rotational difference families that determine seven Steiner system $S(2,6,226)$ is the semidirect product $(\IZ_5\times\IZ_5\times\IZ_3)\rtimes\IZ_3$, endowed with the group operation \[(x,a)\odot(y,b)=(x+y \cdot M^a,a+b),\]
    where $x,y\in \IZ_5\times\IZ_5\times\IZ_3$, $a,b\in\IZ_3$, and
    $M$ is the automorphism of the group $\IZ_5\times\IZ_5\times\IZ_3$, defined by the matrix multiplication
    \[x\cdot M=(x_1, x_2, x_3)\cdot\left(\begin{array}{ccc}4&4&0\\
        1&0&0\\
        0&0&1\end{array}\right)=(4 x_1+x_2,4 x_1,x_3)\]
    where $x=(x_1,x_2,x_3)\in \IZ_5\times\IZ_5\times\IZ_3$.

    \smallskip
    By computer search, we found seven $(226,6,1)$ $1$-rotational difference families on\break $\big((\IZ_5\times\IZ_5\times\IZ_3)\rtimes\IZ_3\big)\cup\{\infty\}$, determining seven Steiner system $S(2,6,226)$.
    Each of these difference families consists of $4$ blocks with trivial stabilizers, $10$ blocks with stabilizers of order $3$, and one block with stabilizer of order $5$.
    So, the total number of blocks equals $4\cdot\frac{225}{1} + 10\cdot\frac{225}{3} + 1\cdot\frac{225}{5} = 1695 = \frac{226(226-1)}{6(6-1)}$.
    \begin{enumerate}
        \item $\{0000, 0001, 1022, 2012, 3012, 4022\}, \{0000, 0010, 1301, 2111, 3411, 4201\}, \newline \{0000, 0011, 0220, 0320, 1101, 4401\}, \{0000, 0012, 0122, 0422, 2210, 3310\}, \newline \{0000, 0101, 1102, 3210, 3311, 4312\}, \{0000, 0201, 1121, 1420, 2202, 3122\}, \newline \{0000, 0401, 1212, 2211, 2310, 4402\}, \{0000, 0301, 2422, 3302, 4120, 4421\}, \newline \{0000, 0102, 0211, 1210, 1312, 4001\}, \{0000, 0202, 0421, 2122, 2420, 3001\}, \newline \{0000, 0121, 0302, 2001, 3120, 3422\}, \{0000, 0311, 0402, 1001, 4212, 4310\}, \newline \{0000, 0111, 1122, 1200, 1311, 2322\}, \{0000, 0221, 2121, 2212, 2400, 4112\}, \newline \{0000, 1400, 2300, 3200, 4100, \infty\};$ \smallskip
        \item $\{0000, 0001, 1022, 2012, 3012, 4022\}, \{0000, 0010, 1402, 2312, 3212, 4102\}, \newline \{0000, 0011, 0110, 0410, 2201, 3301\}, \{0000, 0012, 0202, 0302, 1110, 4410\}, \newline \{0000, 0101, 0312, 1102, 4210, 4311\}, \{0000, 0301, 0422, 2120, 2421, 3302\}, \newline \{0000, 0122, 0201, 2202, 3121, 3420\}, \{0000, 0212, 0401, 1211, 1310, 4402\}, \newline \{0000, 0111, 1122, 1200, 1311, 2322\}, \{0000, 0221, 2121, 2212, 2400, 4112\}, \newline \{0000, 0121, 2410, 2422, 3211, 4302\}, \{0000, 0211, 1421, 3102, 4312, 4320\}, \newline \{0000, 0421, 1202, 2311, 3110, 3122\}, \{0000, 0311, 1212, 1220, 2402, 4121\}, \newline \{0000, 1300, 2100, 3400, 4200, \infty\};$ \smallskip
        \item $\{0000, 0001, 1022, 2012, 3012, 4022\}, \{0000, 0010, 1401, 2311, 3211, 4101\}, \newline \{0000, 0011, 0220, 0320, 1101, 4401\}, \{0000, 0012, 0122, 0422, 2210, 3310\}, \newline \{0000, 0101, 1102, 2110, 2211, 3212\}, \{0000, 0201, 1422, 2202, 4220, 4421\}, \newline \{0000, 0401, 2312, 3311, 3410, 4402\}, \{0000, 0301, 1121, 1320, 3302, 4122\}, \newline \{0000, 0102, 0211, 1210, 1312, 4001\}, \{0000, 0202, 0421, 2122, 2420, 3001\}, \newline \{0000, 0121, 0302, 2001, 3120, 3422\}, \{0000, 0311, 0402, 1001, 4212, 4310\}, \newline \{0000, 0111, 1122, 1200, 1311, 2322\}, \{0000, 0221, 2121, 2212, 2400, 4112\}, \newline \{0000, 1300, 2100, 3400, 4200, \infty\};$ \smallskip
        \item $\{0000, 0001, 1022, 2012, 3012, 4022\}, \{0000, 0010, 1412, 2302, 3202, 4112\}, \newline \{0000, 0011, 0122, 0422, 1212, 4312\}, \{0000, 0012, 0221, 0321, 2401, 3101\}, \newline \{0000, 0101, 0311, 0320, 3112, 3322\}, \{0000, 0110, 0121, 0201, 1112, 1222\}, \newline \{0000, 0120, 1321, 1411, 3212, 3302\}, \{0000, 0211, 0220, 0401, 2222, 2412\}, \newline \{0000, 0102, 1322, 2421, 3111, 4310\}, \{0000, 0202, 1221, 2112, 3120, 4311\}, \newline \{0000, 0402, 1210, 2411, 3121, 4222\}, \{0000, 0302, 1211, 2420, 3412, 4321\}, \newline \{0000, 1002, 1101, 2100, 3102, 3201\}, \{0000, 1202, 1401, 2002, 2201, 4200\}, \newline \{0000, 1100, 2200, 3300, 4400, \infty\};$ \smallskip
        \item $\{0000, 0001, 1022, 2012, 3012, 4022\}, \{0000, 0010, 1411, 2301, 3201, 4111\}, \newline \{0000, 0011, 1101, 1310, 4210, 4401\}, \{0000, 0012, 2102, 2210, 3310, 3402\}, \newline \{0000, 0221, 2121, 2212, 2400, 4112\}, \{0000, 0101, 0110, 0211, 1102, 1212\}, \newline \{0000, 0201, 0220, 0421, 2202, 2422\}, \{0000, 0120, 1301, 1421, 3202, 3322\}, \newline \{0000, 0121, 0301, 0320, 3122, 3302\}, \{0000, 0102, 1322, 2421, 3111, 4310\}, \newline \{0000, 0202, 1221, 2112, 3120, 4311\}, \{0000, 0402, 1210, 2411, 3121, 4222\}, \newline \{0000, 0302, 1211, 2420, 3412, 4321\}, \{0000, 0111, 1122, 1200, 1311, 2322\}, \newline \{0000, 0100, 0200, 0300, 0400, \infty\};$ \smallskip
        \item $\{0000, 0001, 1022, 2012, 3012, 4022\}, \{0000, 0010, 1112, 2202, 3302, 4412\}, \newline \{0000, 0011, 2322, 2412, 3112, 3222\}, \{0000, 0012, 1201, 1421, 4121, 4301\}, \newline \{0000, 0100, 2002, 2402, 3201, 3301\}, \{0000, 0200, 1101, 1401, 4002, 4302\}, \newline \{0000, 0102, 1322, 2421, 3111, 4310\}, \{0000, 0202, 1221, 2112, 3120, 4311\}, \newline \{0000, 0402, 1210, 2411, 3121, 4222\}, \{0000, 0302, 1211, 2420, 3412, 4321\}, \newline \{0000, 0111, 0212, 1122, 4120, 4201\}, \{0000, 0221, 0422, 2212, 3210, 3401\}, \newline \{0000, 0122, 0321, 2101, 2310, 3312\}, \{0000, 0312, 0411, 1301, 1420, 4422\}, \newline \{0000, 1300, 2100, 3400, 4200, \infty\};$ \smallskip
        \item $\{0000, 0001, 1022, 2012, 3012, 4022\}, \{0000, 0010, 0111, 0201, 0301, 0411\}, \newline \{0000, 0011, 1202, 2122, 3422, 4302\}, \{0000, 0012, 1321, 2411, 3111, 4221\}, \newline \{0000, 0102, 3101, 4001, 4100, 4202\}, \{0000, 0202, 1201, 3001, 3200, 3402\}, \newline \{0000, 0121, 1112, 2120, 2211, 3202\}, \{0000, 0211, 1402, 2222, 4210, 4421\}, \newline \{0000, 0421, 2302, 3311, 3420, 4412\}, \{0000, 0311, 1121, 1310, 3322, 4102\}, \newline \{0000, 0122, 1421, 2002, 2410, 4011\}, \{0000, 0212, 2311, 3021, 4002, 4320\}, \newline \{0000, 0422, 1011, 3002, 3110, 4121\}, \{0000, 0312, 1002, 1220, 2021, 3211\}, \newline \{0000, 1100, 2200, 3300, 4400, \infty\}.$
    \end{enumerate}

    The base blocks consist of elements of the group $G=(\IZ_5\times\IZ_5\times\IZ_3)\rtimes\IZ_3$, denoted by quadruples of numbers (without commas and parentheses) or the symbol $\infty$ (which is the fixed point of the action of the group $G$ on $G\cup\{\infty\}$).

    Using \texttt{traces}~\cite{nauty} it was checked that the automorphism groups of all seven $S(2,6,226)$ have order $900$ and are isomorphic to \texttt{SmallGroup(900,92)} from \texttt{GAP}~\cite{GAP4}, that was the group used in the algorithm~\cite{Het} for producing these designs.
    These seven designs exhaust all $1$-rotational Steiner systems $S(2,6,226)$ with automorphism group isomorphic to \texttt{SmallGroup(900,92)}.

    \section{Steiner systems $S(2,6,441)$}\label{sec:441}

    The construction described below is similar to that of Mills~\cite{WHM78} for $S(2,6,111)$.

    By computer search, we found six Steiner systems $S(2,6,441)$, generated by six $(441,6,1)$ difference families on the group $(\IZ_7\rtimes\IZ_3)\times(\IZ_7\rtimes \IZ_3)$, which is the Cartesian square of the group $\IZ_7 \rtimes \IZ_3$, endowed with the group operation $(x,a)\cdot(y,b)=(x+2^a y,a+b)$, where $x,y\in \IZ_7$ and $a,b\in\IZ_3$.
    Each of these difference families consists of $20$ blocks: the first $12$ with trivial stabilizers and the last $8$ with stabilizers of order $3$.
    So, the total number of blocks is $12\cdot\frac{441}{1} + 8\cdot\frac{441}{3} = 6468 = \frac{441(441-1)}{6(6-1)}$.
    \begin{enumerate}
        \item $\{0000, 0001, 1022, 3202, 4202, 6052\}, \{0000, 0010, 0111, 0211, 2001, 5021\}, \newline \{0000, 0011, 0230, 3232, 5131, 6162\}, \{0000, 0042, 0212, 1131, 2150, 4261\}, \newline \{0000, 0110, 2012, 2241, 3101, 5212\}, \{0000, 0160, 2262, 4101, 5062, 5231\}, \newline \{0000, 0020, 1052, 2240, 5250, 6032\}, \{0000, 1040, 1251, 3041, 3231, 5041\}, \newline \{0000, 0030, 1110, 1201, 6120, 6261\}, \{0000, 0112, 1162, 4000, 5061, 6011\}, \newline \{0000, 0101, 3042, 3230, 4032, 4240\}, \{0000, 0102, 1051, 2212, 5262, 6021\}, \newline \{0000, 0012, 0051, 4011, 4030, 4042\}, \{0000, 1130, 1200, 3030, 4230, 5100\}, \newline \{0000, 0022, 0031, 3012, 3021, 3060\}, \{0000, 1030, 4100, 5130, 5200, 6230\}, \newline \{0000, 0122, 1242, 2041, 3250, 5161\}, \{0000, 1141, 2061, 3102, 3252, 5220\}, \newline \{0000, 0152, 2111, 4220, 5031, 6232\}, \{0000, 1151, 2032, 3120, 3212, 5201\};$ \smallskip
        \item $\{0000, 0001, 1062, 3202, 4202, 6012\}, \{0000, 0010, 0111, 0211, 1021, 6001\}, \newline \{0000, 0011, 1212, 2111, 4032, 4142\}, \{0000, 0012, 1142, 2150, 3212, 6131\}, \newline \{0000, 0021, 1112, 4200, 5141, 6140\}, \{0000, 0042, 3061, 3151, 5100, 6211\}, \newline \{0000, 0020, 1161, 2250, 5240, 6151\}, \{0000, 1030, 3031, 4102, 5031, 5152\}, \newline \{0000, 0030, 0150, 0262, 2231, 5231\}, \{0000, 0101, 1251, 2112, 5162, 6221\}, \newline \{0000, 0112, 2032, 3022, 3142, 5030\}, \{0000, 1010, 2031, 3130, 6061, 6150\}, \newline \{0000, 0022, 0031, 4012, 4021, 4060\}, \{0000, 1040, 4100, 5140, 5200, 6240\}, \newline \{0000, 0131, 0222, 3012, 3160, 3221\}, \{0000, 1052, 4131, 5130, 5222, 6261\}, \newline \{0000, 0141, 0252, 4062, 4110, 4251\}, \{0000, 1031, 4101, 5132, 5202, 6230\}, \newline \{0000, 1020, 4111, 5131, 5232, 6252\}, \{0000, 1131, 2060, 3121, 3222, 5212\};$ \smallskip
        \item $\{0000, 0001, 1042, 3202, 4202, 6032\}, \{0000, 0010, 0111, 0211, 3021, 4001\}, \newline \{0000, 0011, 2001, 2040, 2240, 6200\}, \{0000, 0012, 0150, 1100, 4050, 4052\}, \newline \{0000, 0020, 0121, 0210, 1142, 6142\}, \{0000, 0022, 1162, 2110, 4111, 6251\}, \newline \{0000, 0122, 3222, 4021, 4130, 4220\}, \{0000, 0152, 3051, 3140, 3250, 4252\}, \newline \{0000, 0041, 1232, 3102, 5131, 6130\}, \{0000, 0030, 2022, 2102, 5032, 5152\}, \newline \{0000, 0112, 2062, 3032, 3112, 5040\}, \{0000, 0101, 1061, 1110, 6011, 6160\}, \newline \{0000, 1010, 4141, 5151, 5252, 6262\}, \{0000, 1060, 4131, 5121, 5222, 6212\}, \newline \{0000, 1030, 4112, 5142, 5251, 6211\}, \{0000, 1112, 1231, 3060, 4221, 5122\}, \newline \{0000, 1111, 2012, 3150, 3232, 5261\}, \{0000, 1161, 1250, 3242, 5041, 6122\}, \newline \{0000, 1152, 2031, 4232, 6111, 6220\}, \{0000, 1141, 2021, 3132, 3252, 5250\};$ \smallskip
        \item $\{0000, 0001, 1022, 3201, 4201, 6052\}, \{0000, 0010, 0011, 0211, 2101, 5121\}, \newline \{0000, 0012, 3242, 5240, 6041, 6160\}, \{0000, 0160, 2252, 3102, 4261, 5260\}, \newline \{0000, 0110, 2210, 3211, 4102, 5222\}, \{0000, 0021, 1002, 1131, 2251, 4240\}, \newline \{0000, 0030, 3012, 3210, 4042, 4220\}, \{0000, 0112, 1130, 3060, 5130, 6122\}, \newline \{0000, 0101, 2051, 2120, 5021, 5150\}, \{0000, 0102, 1051, 2212, 5262, 6021\}, \newline \{0000, 1020, 1232, 3252, 5231, 6211\}, \{0000, 1112, 2060, 3122, 3262, 5242\}, \newline \{0000, 0020, 0111, 0131, 0232, 0252\}, \{0000, 1101, 2000, 3101, 3202, 5202\}, \newline \{0000, 0022, 0031, 4012, 4021, 4060\}, \{0000, 1040, 4100, 5140, 5200, 6240\}, \newline \{0000, 0130, 2111, 4241, 5062, 6232\}, \{0000, 0140, 1242, 2012, 3231, 5161\}, \newline \{0000, 1250, 3142, 4001, 4262, 6121\}, \{0000, 1151, 3001, 3212, 4132, 6220\};$ \smallskip
        \item $\{0000, 0001, 1062, 3200, 4200, 6012\}, \{0000, 0010, 0011, 0111, 1221, 6201\}, \newline \{0000, 0012, 1141, 2240, 5120, 5231\}, \{0000, 0021, 2101, 2212, 5251, 6102\}, \newline \{0000, 0121, 0261, 2152, 5100, 6222\}, \{0000, 0130, 0212, 1121, 2012, 5061\}, \newline \{0000, 0020, 1160, 2032, 5052, 6130\}, \{0000, 0131, 1120, 2040, 3120, 4151\}, \newline \{0000, 0030, 3032, 3212, 4022, 4242\}, \{0000, 1010, 2052, 3140, 6062, 6140\}, \newline \{0000, 0101, 1142, 3211, 4261, 6132\}, \{0000, 0102, 2051, 2222, 5021, 5252\}, \newline \{0000, 0022, 0031, 4012, 4021, 4060\}, \{0000, 1040, 4100, 5140, 5200, 6240\}, \newline \{0000, 0110, 2121, 4231, 5002, 6262\}, \{0000, 0160, 1212, 2002, 3241, 5151\}, \newline \{0000, 0120, 1002, 4252, 5261, 6141\}, \{0000, 0150, 1131, 2211, 3222, 6002\}, \newline \{0000, 1020, 4132, 5152, 5211, 6231\}, \{0000, 1122, 2060, 3112, 3231, 5221\};$ \smallskip
        \item $\{0000, 0001, 1062, 3202, 4202, 6012\}, \{0000, 0010, 0111, 0211, 1021, 6001\}, \newline \{0000, 0011, 1160, 2201, 3012, 4220\}, \{0000, 0042, 3222, 4011, 5210, 6152\}, \newline \{0000, 0122, 1011, 5042, 6100, 6230\}, \{0000, 0112, 3232, 4151, 4230, 5130\}, \newline \{0000, 0012, 2221, 4261, 5012, 5040\}, \{0000, 0020, 2120, 2232, 5100, 5252\}, \newline \{0000, 0030, 3001, 3160, 4061, 4140\}, \{0000, 0130, 1031, 1140, 3041, 4000\}, \newline \{0000, 0101, 2042, 3220, 4250, 5032\}, \{0000, 0102, 2151, 2220, 5121, 5250\}, \newline \{0000, 0022, 0031, 2012, 2021, 2060\}, \{0000, 1020, 4100, 5120, 5200, 6220\}, \newline \{0000, 0120, 2202, 3101, 4022, 6221\}, \{0000, 0150, 1251, 3052, 4101, 5202\}, \newline \{0000, 0131, 0222, 5001, 5162, 5240\}, \{0000, 0141, 0252, 2001, 2112, 2230\}, \newline \{0000, 1040, 4121, 5161, 5262, 6232\}, \{0000, 1151, 1252, 3010, 4262, 5141\}.$
    \end{enumerate}

    Writing the generating blocks, we denoted elements of the group $(\IZ_7\rtimes\IZ_3)\times(\IZ_7\rtimes\IZ_3)$ by quadruples of numbers (without commas and parentheses).

    Using \texttt{traces}~\cite{nauty} it was checked that the automorphism groups of all six $S(2,6,441)$ have order $1764$ and are isomorphic to \texttt{SmallGroup(1764,133)} from \texttt{GAP}~\cite{GAP4}, that was the group used in the algorithm~\cite{Het} for producing these designs.
    Comparing to the previous case, the search was not exhaustive, lasted one hour without even finishing one initial search state out of millions.
    This means that this group is potentially profilic for generating many $S(2,6,441)$.

\end{document}